\documentclass[10pt]{article}
\usepackage{amssymb}
\begin{document}
\title{The 2-color Rado Number of $x_1+x_2+\cdots +x_n=y_1+y_2+\cdots +y_k$}
\author{Dan Saracino\\Colgate University}
\date{}
\maketitle
\begin{abstract}  In 1982, Beutelspacher and Brestovansky determined the 2-color Rado number of the equation $$x_1+x_2+\cdots +x_{m-1}=x_m$$ for all $m\geq 3.$ Here we extend their result by determining the 2-color Rado number of the equation $$x_1+x_2+\cdots +x_n=y_1+y_2+\cdots +y_k$$ for all $n\geq 2$ and $k\geq 2.$ As a consequence, we determine the 2-color Rado number of 
$$x_1+x_2+\cdots +x_n=a_1y_1+\cdots +a_{\ell}y_{\ell}$$ in all cases where $n\geq 2$ and $n\geq a_1+\cdots +a_{\ell},$ and in most cases where $n\geq 2$ and $2n\geq a_1+\cdots +a_{\ell}.$
\end{abstract}
\vspace{.25in}

\noindent \textbf{1. Introduction}

\vspace{.25in}

A special case of the work of Richard Rado [\textbf{5}] is that for all positive integers $n$ and $k$ such that $n+k\geq 3$, and  all positive integers $a_1,\ldots,a_n$ and $b_1,\ldots,b_k,$ there exists a smallest positive integer $r$ with the following property:  for every coloring of the elements of the set $[r]=\{1,\ldots,r\}$ with two colors, there  exists a solution of the equation $$a_1x_1+a_2x_2+\cdots +a_nx_n=b_1y_1+b_2y_2+\cdots +b_ky_k$$ using elements of $[r]$ that are all colored the same. (Such a solution is called \emph{monochromatic}.) The integer $r$ is called the \emph{2-color Rado number} of the equation.

In recent years there has been a considerable amount of work aimed at determining the Rado numbers of specific equations.  One of the earliest results in this direction
appeared in a 1982 paper of Beutelspacher and Brestovansky [\textbf{1}], where it was proved that for every $m\geq 3$, the 2-color Rado number of $$x_1+x_2+\cdots +x_{m-1}=x_m$$ is $m^2-m-1.$  In 2008 Guo and Sun [\textbf{2}] generalized this result by proving that, for all positive integers $a_1,\ldots,a_{m-1},$ the 2-color Rado number of the equation $$a_1x_1+a_2x_2+\cdots+a_{m-1}x_{m-1}=x_m$$  is $aw^2+w-a,$ where $a=\textrm{min}\{a_1,\dots, a_{m-1}\}$ and $w=a_1+\cdots +a_{m-1}.$
In the same year, Schaal and Vestal $[\textbf{8}]$ dealt with the equation $$x_1+x_2+\cdots +x_{m-1}=2x_m.$$ They proved, in particular, that for every $m\geq 6,$ the 2-color Rado number is $\lceil\frac{m-1}{2}\lceil\frac{m-1}{2}\rceil\rceil.$  Building on the work of Schaal and Vestal, we investigated the equation

\begin{equation}x_1+x_2+\cdots +x_{m-1}=ax_m, \end{equation}
 for $a\geq 3,$ in $[\textbf{6}]$ and $[\textbf{7}].$

\vspace{.15in}

\noindent\textbf{Notation.}  We will denote $\lceil\frac{m-1}{a}\lceil\frac{m-1}{a}\rceil\rceil$ by $C(m,a),$ and we will denote the 2-color Rado number of equation (1) by $R_2(m,a).$

\vspace{.15in}

\noindent\textbf{Fact 1} ($[\textbf{7}]$, Theorem 3).  Suppose $a\geq 3.$ If $3|a$ then $R_2(m,a)=C(m,a)$ when $m\geq 2a+1$ but $R_2(2a,a)=5$ and $C(2a,a)=4.$ If $3\nmid a$ then $R_2(m,a)=C(m,a)$  when $m\geq 2a+2$ but $R_2(2a+1,a)=5$ and \linebreak$C(2a+1,a)=4.$

\vspace{.15in}

We note that, by the results of $[\textbf{8}]$, the statements in Fact 1 remain valid when $a=2.$

Results have been obtained for a number of other variations of the equation $x_1+\cdots +x_{m-1}=x_m,$   most of which have had the property that one side of the equation involves only one variable.  Our first purpose here is to determine the 2-color Rado number of the equation 
\begin{equation} x_1+x_2+\cdots +x_n=y_1+y_2+\cdots +y_k, \end{equation}
for all $n\geq 2$ and $k\geq 2.$  

\vspace{.15in}

\noindent\textbf{Notation.}  We denote the 2-color Rado number of equation (2) by $r_2(n,k).$

\vspace{.15in} To determine $r_2(n,k)$ for all $n$ and $k$ it clearly suffices to consider the case $n\geq k.$  We will deal with this case by relating equation (2) to the equation \begin{equation}x_1+\cdots +x_n=ky\end{equation} and using results from $[\textbf{7}]$ (and $[\textbf{1}]$ and $[\textbf{8}]$, for the cases $k=1,2$). 

\vspace{.15in}

\noindent\textbf{Theorem 1.}  If $n\geq 2$ and $n\geq k,$ then $r_2(n,k)=R_2(n+1,k).$

\vspace{.15in}

The relevant values of $R_2(n+1,k)$  are determined by $[\textbf{1}]$, $[\textbf{8}]$, Fact 1 and the following additional information from $[\textbf{7}]$. 

\vspace{.15in}

\noindent\textbf{Fact 2}  ($[\textbf{7}]$, Theorem 2).  Suppose $a+1\leq m\leq 2a+1.$ Then $R_2(m,a)=1$ iff  $m=a+1.$  If $a+2\leq m\leq 2a+1,$ then $R_2(m,a)\in \{3,4,5\},$ and we have:

\begin{itemize}

\item[] $R_2(m,a)=3$  iff $m\leq \frac{3a}{2}+1$ and $a\equiv m-1$ (mod 2).

\item[] $R_2(m,a)=4$ iff either:

\begin{itemize}
\item[(i)] $m\leq \frac{3a}{2}+1$ and $a\not \equiv m-1$ (mod 2), or

\item[(ii)] $m> \frac{3a}{2}+1$ and $a\equiv m-1$  (mod 3).

\end{itemize}

\item[] $R_2(m,a)=5$ iff $m> \frac{3a}{2}+1$ and $a\not \equiv m-1$  (mod 3).

\end{itemize}

\vspace{.15in}

The result of Theorem 1 is not valid, in general, if $n< k.$  We provide the following additional information.

\vspace{.15in}

\noindent\textbf{Theorem 2.}  If $n\geq 2$ and $2n\geq  k>  n,$ then $r_2(n,k)=R_2(n+1,k)$ except in the following cases:

\begin{itemize}

\item[]  $r_2(2,3)=4$ while $R_2(3,3)=9,$

\item[]  $r_2(2,4)=5$ while $R_2(3,4)=10,$

\item[] $r_2(3,5)=5$ while $R_2(4,5)=9,$ and

\item[] if $10\leq k \leq 14$ then $r_2(k-5,k)=5$ while $R_2(k-4,k)=6.$

\end{itemize}

\vspace{.15in}

The proof of Theorem 2 relies on the following two results from $[\textbf{7}]$.

\vspace{.15in}

\noindent\textbf{Fact 3}  ($[\textbf{7}]$, Theorem 4). If  $\frac{2a}{3}+1\leq m\leq a,$ then:

\begin{itemize} \item[] for $a=3$ we have $R_2(a,a)=9,$ and \item[] for  $a\geq 4$ we have

\begin{itemize} \item[] $R_2(m,a)=3$ if $a\equiv m-1$ (mod 2) and \item[] $R_2(m,a)=4$ if $a \not \equiv m-1$  (mod 2).  

\end{itemize}
\end{itemize}

\vspace{.15in}

\noindent\textbf{Fact 4} ($[\textbf{7}],$ Theorem 5). If $\frac{a}{2}+1\leq m < \frac{2a}{3}+1$ (so $a\geq 4)$  then:

\begin{itemize}
\item[] for $a\equiv m-1$  (mod 3) we have $R_2(m,a)=4,$ and

\item[] for $a\not \equiv m-1$ (mod 3) we have $R_2(m,a)=5$ \emph{except} that
\begin{itemize}
\item[]  $R_2(3,4)=10$ and $R_2(4,5)=9,$ and

\item[] $R_2(m,a)=6$ if $10\leq a\leq 14$ and $m=a-4.$ 

\end{itemize}

\end{itemize}

\vspace{.15in}

We will show that Theorems 1 and 2 have the following consequence.

\vspace{.15in}

\noindent\textbf{Theorem 3.} Let $n\geq 2,$ let  $a_1,\ldots, a_{\ell}$ be positive integers, and let $A=a_1+\cdots +a_{\ell}.$ Then if $n\geq A,$ the 2-color Rado number of $$x_1+x_2+\cdots +x_n=a_1y_1+\cdots +a_{\ell}y_{\ell}$$ is $R_2(n+1,A).$  If $2n\geq A> n,$ the same conclusion holds provided that the pair $(n,A)$ is none of $(2,3), (2,4), (3,5), (5,10), (6,11), (7,12), (8,13), (9,14). $

\vspace{.15in}

For the case $A> 2n$ we have the following.  \vspace{.15in}

\noindent\textbf{Theorem 4.} If $n\geq 2, \ A=a_1+\cdots +a_{\ell}$ and $A>2n,$ then the 2-color Rado number of $$x_1+x_2+\cdots +x_n=a_1y_1+\cdots +a_{\ell}y_{\ell}$$ is at least $\lceil\frac{A}{n}\lceil\frac{A}{n}\rceil\rceil.$

 \vspace{.25in}
 
 \noindent\textbf{2.  Proofs of the Theorems}
 
 \vspace{.25in}
 
 \noindent\textbf{Lemma.}  For any $n\geq k,$ we have $$C(n+1,k)\leq r_2(n,k)\leq R_2(n+1,k).$$
 
 \vspace{.15in}
 
 \noindent\emph{Proof.} The second inequality is clear, since any monochromatic solution of equation (3) provides a monochromatic solution of equation (2).
 
 To prove the first inequality it suffices to consider $n> k$ and exhibit a 2-coloring of $[C(n+1,k)-1]$ that yields no monochromatic solution of equation (2). Let the elements of $[\lceil\frac{n}{k}\rceil-1]$ be colored red and let the remaining elements of $[C(n+1,k)-1]$ be colored blue. For any red solution of equation (2), the left side has total value at least $n$ and the right side has total value at most $k(\lceil\frac{n}{k}\rceil-1) < n,$ so there is no red solution. For any blue solution the left side has total value at least $n\lceil\frac{n}{k}\rceil$ and the right side has total value at most $k(\lceil\frac{n}{k}\lceil\frac{n}{k}\rceil\rceil-1) < n\lceil\frac{n}{k}\rceil$, so there is no blue solution.  $\Box$
 
 \vspace{.15in}
 
 \noindent\emph{Proof of Theorem 1}.  Theorem 1 is obviously true when $k=1,$ so we assume $k\geq 2.$ 
 
 If $3|k$ and $n\geq 2k$ or if $3\nmid k$ and $n\geq 2k+1$ then by Fact 1 and $[\textbf{8}]$ we have $R_2(n+1,k)=C(n+1,k),$ so by the Lemma we have $r_2(n,k)=R_2(n+1,k).$  To complete the proof, it suffices to consider the cases where $k\leq n\leq 2k.$
 
 Note that $r_2(n,k)=1$ iff $n=k$ iff $R_2(n+1,k)=1,$ so we can suppose that $k+1\leq n\leq 2k.$ Then by coloring 1 and 2 differently, we see that $r_2(n,k)$ cannot be 2, so  $r_2(n,k)\geq 3.$  We have $R_2(n+1,k)\in \{3,4,5\}$ by Fact 2.   We now consider the mutually exclusive cases indicated in Fact 2.
 
 First suppose that $n\leq \frac{3k}{2}$ and $k\equiv n$ (mod 2). Then by Fact 2 we have $R_2(n+1,k)=3.$ By the Lemma, $r_2(n,k)$ can only be 3.
 
 Next suppose that $n\leq \frac{3k}{2}$ and $k\not \equiv n$ (mod 2). Then $R_2(n+1,k)=4$ by Fact 2.  By the Lemma, $r_2(n,k)=3$ or $4$.  If we 2-color $[3]$  by coloring 1 and 3 red and 2 blue, then since $k\not \equiv n$ (mod 2) there is no red solution of equation (2), and since $k\neq n$ there is no blue solution. So $r_2(n,k)=4.$
 
Now suppose that $n> \frac{3k}{2}$ and $k\equiv n$ (mod 3).  Then by Fact 2 we have $R_2(n+1,k)=4,$ so again $r_2(n,k)=3$ or $4$.  Since $n>\frac{3k}{2},$ we have $C(n+1,k)=\lceil\frac{n}{k}\lceil\frac{n}{k}\rceil\rceil \geq \lceil\frac{n}{k}\cdot 2\rceil$ and $\frac{n}{k}\cdot 2>3,$ so it follows from the Lemma that $r_2(n,k)\geq 4,$ and therefore $r_2(n,k)=4.$

 Finally, if $n>\frac{3k}{2}$ and $k\not \equiv n$ (mod 3), then $R_2(n+1,k)=5$ by Fact 2.  As in the preceding paragraph, we have $r_2(n,k)\geq 4.$  If we 2-color $[4]$ by coloring 1 and 4 red and 2 and 3 blue, then for any blue solution of equation (2) the left side has total value at least $2n$  and the right side has total value at most $3k$. Since $2n> 3k,$ there is no blue solution.  For any red solution the left side of the equation has total value congruent to $n$ (mod 3) and the right side has total value congruent to $k$ (mod 3). Since $k\not \equiv n$ (mod 3) there can be no red solution. So $r_2(n,k)=5=R_2(n+1,k).$  $\Box$
 
 \vspace{.15in}
 
 \noindent\emph{Proof of Theorem 2.}  Since $k> n,$ we have $r_2(n,k)=r_2(k,n)=R_2(k+1,n),$ by Theorem 1.  To evaluate $R_2(k+1,n)$ we will use Fact 2.  Note that we can do so since the condition $a+2\leq m\leq 2a+1$ of Fact 2, with $a=n$ and $m=k+1,$ becomes $n+2\leq k+1\leq 2n+1,$ which holds since $2n\geq k> n.$
 
 \vspace{.1in}
 \noindent\emph{Case 1: $k> n\geq  \frac{2k}{3}.$}
 
 \vspace{.1in}
 
 In this case, if $k=3$ then $n=2$.  We have $r_2(2,3)=R_2(3+1,2)=4$ by $[\textbf{8}]$, while $R_2(n+1,k)=R_2(3,3)=9$ by Fact 3.
 
 We claim that if $k\geq 4$ then $r_2(n,k)=R_2(n+1,k),$ i.e., $R_2(k+1,n)=R_2(n+1,k).$
 Note that since $k\leq \frac{3n}{2},$ Fact 2 (with $a=n$ and $m=k+1$) yields $R_2(k+1,n)=3$ if $n\equiv k$  (mod 2) and $R_2(k+1,n)=4$ if $n\not \equiv k$ (mod 2).  To determine $R_2(n+1,k)$ we can use Fact 3 (with $a=k$ and $m=n+1$), since $\frac{2k}{3}+1\leq n+1\leq k$.  We find that $R_2(n+1,k)=3$ if $k\equiv n$ (mod 2) and $R_2(n+1,k)=4$ if $k\not \equiv n$ (mod 2).  This concludes the proof in Case 1.
 
 \vspace{.1in}
 
 \noindent\emph{Case 2: $\frac{2k}{3}> n\geq \frac{k}{2}.$}
 
 \vspace{.1in}
 
 Using Fact 2 with $a=n$ and $m=k+1$, we note that since $k>\frac{3n}{2}$ we have $R_2(k+1,n)=4$ if $n\equiv k$ (mod 3) and $R_2(k+1,n)=5$ if $n \not \equiv k$ (mod 3).  Using Fact 4 with $a=k$ and $m=n+1$ (which is legitimate since $\frac{k}{2}+1\leq n+1< \frac{2k}{3}+1$) we find that $R_2(n+1,k)=4$ if $k\equiv n$ (mod 3) and $R_2(n+1,k)=5$ if $k\not \equiv n$ (mod 3) (and therefore $r_2(n,k)=R_2(n+1,k)$) \emph{unless} the pair $(n+1,k)$ is $(3,4)$, $(4,5),$ or $(k-4,k)$ for some $10\leq k\leq 14,$ in which case $R_2(n+1,k)$ is 10,  9, or 6, respectively.  $\Box$

 \vspace{.15in}
 
 \noindent\emph{Proof of Theorem 3.}  The 2-color Rado number of
 
 \begin{equation} x_1+x_2+\cdots +x_n=a_1y_1+\cdots +a_{\ell}y_{\ell} \end{equation} is at least that of
 \begin{equation} x_1+x_2+\cdots +x_n=y_1+\cdots +y_A \end{equation} and at most that of
 \begin{equation} x_1+x_2+\cdots +x_n=Ay, \end{equation}
 since any monochromatic solution of equation (4) yields a monochromatic solution of equation (5) and any monochromatic solution of equation (6) yields a monochromatic solution of equation (4).   The 2-color Rado numbers of equations (5) and (6) are the same by Theorem 1 if $n\geq A,$  and they are the same by Theorem 2, with the stated exceptions, when $2n\geq A > n.$ Thus the 2-color Rado number of equation (4) is that of equation (6), namely $R_2(n+1,A),$ if $(n,A)$ is not one of the indicated exceptional pairs.   $\Box$
 
 \vspace{.15in}
 
 We have not ruled out the possibility that the 2-color Rado number for some instances  of equation (4) with $\ell>1$ is $R_2(n+1,A)$ even when $(n,A)$ is one of the exceptional pairs.
 
 \vspace{.15in}
 
 \noindent\emph{Proof of Theorem 4.}  As above, the 2-color Rado number of equation (4) is at least that of equation (5), which is $R_2(A+1,n)$ by Theorem 1.  Since $A> 2n,$ we have $A+1\geq 2n+2,$ so $R_2(A+1,n)=\lceil\frac{A}{n}\lceil\frac{A}{n}\rceil\rceil$ by Fact 1.  $\Box$

 \vspace{.25in}

\centerline{ \textbf{References}}
\vspace{.1in}

\noindent 1.  A. Beutelspacher and W. Brestovansky, Generalized Schur numbers, \emph{Lecture Notes in Mathematics}, Springer-Verlag, Berlin, \textbf{969} (1982), 30-38

\vspace{.05in}

\noindent 2. S. Guo and Z-W. Sun, Determination of the 2-color Rado
number for $a_1x_1+ \cdots +a_mx_m=x_0$,  \emph{J. Combin. Theory Ser.A}, \textbf{115} (2008), 345-353.

\vspace{.05in}

\noindent 3.  H. Harborth and S. Maasberg,  All two-color Rado numbers for $a(x+y)=bz,$ \emph{Discrete Math.} \textbf{197-198} (1999), 397-407.

\vspace{.05in}

\noindent 4. B. Hopkins and D. Schaal, On Rado numbers for
$\displaystyle \Sigma_{i=1}^{m-1} a_ix_i=x_m$, \emph{Adv. Applied
Math.} \textbf{35} (2005), 433-441.
\vspace{.05in}

\noindent 5. R. Rado, Studien zur Kombinatorik, \emph{Mathematische
Zeitschrift} \textbf{36} (1933), 424-448.
\vspace{.05in}

\noindent 6. D. Saracino, The 2-color Rado Number of $x_1+x_2+\cdots +x_{m-1}=ax_m,$  \emph{Ars Combinatoria} \textbf{113} (2014), 81-95.
\vspace{.05in}

\noindent 7. D.Saracino, The 2-color Rado number of $x_1+x_2+\cdots +x_{m-1}=ax_m, II$ arXiv://math.CO/1306.0775
\vspace{.05in}

\noindent 8. D. Schaal and D. Vestal, Rado numbers for $x_1+x_2+\cdots +x_{m-1}=2x_m$, \emph{Congressus Numerantium}  \textbf{191} (2008), 105-116.

\end{document}